\documentclass[a4paper,reqno, 11pt]{amsart}

\usepackage[margin=3cm]{geometry}

\makeatletter
\let\old@setaddresses\@setaddresses
\def\@setaddresses{\bigskip\bgroup\parindent 0pt\let\scshape\relax\old@setaddresses\egroup}
\makeatother

\usepackage{amssymb} %Box
\usepackage{amsmath} %Operator
\usepackage{amsthm} %proof
\usepackage[unicode=true]{hyperref}
\hypersetup{
    colorlinks,
    linkcolor={red!50!black},
    citecolor={blue!50!black},
    urlcolor={blue!80!black}
}

\usepackage{color}

\usepackage{todonotes}
\usepackage{bbm}
\usepackage{paralist}
\usepackage[capitalise]{cleveref}

\newtheorem{theorem}{Theorem}[section]

\theoremstyle{remark}
\newtheorem*{acknowledgement}{Acknowledgments}

\makeatletter
\newcommand{\vast}{\bBigg@{4}}
\newcommand{\Vast}{\bBigg@{5}}
\makeatother
\definecolor{bulgarianrose}{rgb}{0.28, 0.02, 0.03}
\definecolor{gray}{rgb}{0.5, 0.5, 0.5}

\makeatletter
\def\namedlabel#1#2{\begingroup
    #2%
    \def\@currentlabel{#2}%
    \phantomsection\label{#1}\endgroup
}
\usepackage[normalem]{ulem}
\usepackage{dsfont}
\usepackage{accents}
\usepackage{verbatim}
\usepackage{ulem}
\usepackage{pgf,tikz,pgfplots}
\pgfplotsset{compat=1.16}%leave out before submission
\usepackage[all]{xy} 

\usepackage{stackengine}
\stackMath
\newcommand\tsup[2][2]{%
 \def\useanchorwidth{T}%
  \ifnum#1>1%
    \stackon[-.5pt]{\tsup[\numexpr#1-1\relax]{#2}}{\scriptscriptstyle\sim}%
  \else%
    \stackon[.5pt]{#2}{\scriptscriptstyle\sim}%
  \fi%
}

\usepackage{subfig}

\usepackage{enumerate}

\title{Local certification of graphs on surfaces}

\author[L.~Esperet]{Louis Esperet}
\address[L.~Esperet]{Laboratoire G-SCOP (CNRS, Univ.\ Grenoble Alpes), Grenoble, France}
\email{louis.esperet@grenoble-inp.fr}

\author[B.~Lévêque]{Benjamin Lévêque}
\address[B. Lévêque]{Laboratoire G-SCOP (CNRS, Univ.\ Grenoble Alpes), Grenoble, France}
\email{benjamin.leveque@grenoble-inp.fr}

\thanks{The authors are partially supported by the French ANR Projects GATO (ANR-16-CE40-0009-01), GrR (ANR-18-CE40-0032), and by LabEx PERSYVAL-lab (ANR-11-LABX-0025). }

\begin{document}

\maketitle
 
\begin{abstract}
A proof labelling scheme for a graph class $\mathcal{C}$ is an assignment of certificates to the vertices of any graph in the class $\mathcal{C}$, such that upon reading its certificate and the certificates of its neighbors, every vertex from a graph $G\in \mathcal{C}$ accepts the instance, while if $G\not\in \mathcal{C}$, for every possible assignment of certificates, at least one vertex rejects the instance.  It was proved recently that for any fixed surface $\Sigma$, the class of graphs embeddable in $\Sigma$ has a proof labelling scheme in which each vertex of an $n$-vertex graph receives a certificate of at most $O(\log n)$ bits. The proof is quite long and intricate and heavily relies on an earlier result for planar graphs. Here we give a very short proof for any surface. The main idea is to encode a rotation system locally, together with a spanning tree supporting the local computation of the genus via Euler's formula.\\

\end{abstract}

%\hspace{1em}Keywords: Bla

%\hspace{1em}MSC Class: Bla

\section{Introduction}

The goal of local certification is to verify that a network, represented by a connected graph $G$ in which each vertex has a unique identifier, satisfies some given property. The constraint is that each node of the network has a local view of the network (its neighborhood) and has to make its decision based only on this local view. If the graph satisfies the property, we want all vertices to accept the instance, while if the graph does not satisfy the property, at least one vertex has to reject the instance. This is a significant restriction and it only allows the verification of local properties (related to the degrees, for instance), so each vertex is given in addition some small certificate, and each vertex can now base its decision on its certificate and the certificates of its neighbors. For any property $\mathcal{P}$, the goal is to produce a protocol to certify $\mathcal{P}$ locally while using certificates of minimal size. Such a protocol is called a \emph{one-round proof labelling scheme with complexity $f(n)$}, where $f(n)$ is the maximum number of bits in the certificate of a vertex in an $n$-vertex graph satisfying $\mathcal{P}$ (a  formal definition of proof labelling schemes  will be given in Section~\ref{sec:pls}). Proof labelling schemes are a natural component of self-stabilizing algorithms, and are a particular form of distributed interactive protocols (with a single interaction). More broadly, proof labelling schemes with compact certificates (of logarithmic or polylogarithmic size) can be seen as a distributed version of the class \textsf{NP}, for which certificates of polynomial size exist (and can be checked in polynomial time in the centralized setting).

\medskip

The Euler genus of a surface $\Sigma$ is denoted by $\mathbf{eg}(\Sigma)$ (see Section~\ref{sec:surfaces} for more on surfaces and graph embeddings).
In particular, the orientable surfaces of Euler genus 0, 2 and 4 are respectively 
the sphere (or equivalently the plane), the torus and the double torus.
The non-orientable surfaces of Euler genus 1 and 2 are respectively the projective plane and the Klein bottle.

\medskip

Motivated by recent work on distributed interactive protocols in classes with linear time recognition algorithms~\cite{NPY20}, it was recently proved that graph planarity has a one-round proof labelling scheme with complexity $O(\log n)$~\cite{planar}, and that this complexity is the best possible. More recently, the same authors built upon their previous work to extend their result to graphs embeddable on any fixed surface~\cite{genus}.

\begin{theorem}[\cite{genus}]\label{thm:main}
For any (orientable or non-orientable) surface $\Sigma$, the class of graphs that are embeddable on $\Sigma$
has a one-round proof labelling scheme with complexity at most $O(\sqrt{\mathbf{eg}(\Sigma)}\cdot \log n)$.
\end{theorem}

%\begin{theorem}[\cite{genus}]\label{thm:main}
%For any integer $g\ge 0$, the classes $\mathcal{C}_g$, %$\mathcal{O}_g$, and $\mathcal{N}_g$ have a one-round proof %labelling scheme with complexity at most $O(\sqrt{g}\cdot \log %n)$.
%\end{theorem}

The proof of the planar case (i.e., the case $\mathbf{eg}(\Sigma)=0$) in~\cite{planar}  and its extension to general surfaces~\cite{genus} are fairly intricate, with the two papers totaling 65 pages. The proof of the planar case~\cite{planar} reduces the problem to graphs that are closer and closer to trees (for which compact proof labelling schemes are known), while the proof for general surfaces~\cite{genus} works by carefully cutting the surface along non-contractible cycles, thus reducing the problem to planar graphs.

In this short note, we give a simple and direct proof of Theorem~\ref{thm:main}, based on rotation systems together with a distributed computation of the Euler genus using Euler's formula along a rooted spanning tree. We believe that our simplified approach is an important step towards an extension of this work to more general classes, such as minor-closed classes.
In addition, we want to emphasize that surfaces are central in the study of distributed algorithms in planar graphs, as these graphs are locally indistinguishable from graphs on surfaces (see for instance~\cite{col} for applications of this observation to obtain lower bounds on distributed coloring of planar graphs). 

\subsubsection*{Related work} A reviewer pointed out an interesting article on a related topic, by Benjamini and Lov\'asz~\cite{BL}, where it is proved that the local observation of some random process in an embedded graph allows to determine the genus of the embedding. The setting is quite different, as the random process is not distributed and the values taken by the edges are real numbers, but it might be the case that some ideas developed there can be useful in the field of local certification.

\subsubsection*{Organization of the paper} The formal definition of proof labelling scheme is given in Section~\ref{sec:pls}, and the basic terminology of graphs on surfaces is given in Section~\ref{sec:surfaces}, along with a description of rotation systems and the Heffter-Edmonds-Ringel rotation principle. For the description of the certificates, we found it convenient to first present the orientable case (Section~\ref{sec:o}), which is slightly simpler, and then explain the small modifications we have to perform in the non-orientable case (Section~\ref{sec:no}). We could have presented everything in the latter setting, which is more general, but we believe it would have been harder to follow. We conclude with some open problems in Section~\ref{sec:ccl}.

\section{Proof labelling schemes}\label{sec:pls}

All graphs in this paper are undirected, simple, and connected.
A \emph{one-round proof labelling scheme} for a graph class $\mathcal{F}$ is a prover-verifier pair $(\mathcal{P},\mathcal{V})$, with the following properties. For any integer $n\ge 1$ and any $n$-vertex graph $G\in \mathcal{F}$, whose vertices are assigned distinct identifiers $(\text{id}(v))_{v\in V(G)}$ from $\{1,\ldots,\text{poly}(n)\}$, the \emph{prover} $\mathcal{P}$ assigns to each vertex $v \in V(G)$ a certificate $c_G(v)\in \{0,1\}^*$ (that might depend on the vertex identifiers). The \emph{verifier} $\mathcal{V}$ satisfies the following properties for any graph $G$:

\medskip

\noindent {\bf One-round:} Each vertex $v\in V(G)$ collects the identifiers and certificates of its neighbors (\emph{one-round}). Let $\mathcal{C}_G(v)=(\text{id}(u),c_G(u))_{uv\in E(G)}$.

\smallskip

\noindent {\bf Completeness:} If $G\in \mathcal{F}$, then for any vertex $v\in V(G)$, \[\mathcal{V}\big(\text{id}(v),c_G(v),\mathcal{C}_G(v)\big)=1.\]

\smallskip

\noindent {\bf Soundness:}  If $G\not\in \mathcal{F}$, then for every possible choice of certificates $(c'_G(v))_{v\in V(G)}$ and distinct identifiers $(\text{id}(v))_{v\in V(G)}$, there exists a vertex $v\in V(G)$ such that  \[\mathcal{V}\big(\text{id}(v),c'_G(v),\mathcal{C}'_G(v)\big)=0,\] where $\mathcal{C}'_G(v)=(\text{id}(u),c'_G(u))_{uv\in E(G)}$.

\medskip

In other words, upon reading the identifiers and certificates of its neighbors and itself, each vertex of a graph $G\in \mathcal{F}$ \emph{accepts} the instance, while if $G\not\in \mathcal{F}$, for every possible choice of identifiers and certificates, at least one vertex  \emph{rejects} the instance. 

\medskip

The \emph{complexity} of the labelling scheme is the maximum size $|c_G(v)|$ of a certificate in an $n$-vertex graph of $\mathcal{F}$. If we say that the complexity is $O(f(n))$, for some function $f$, the $O(\cdot)$ notation refers to $n\to \infty$. The definition above assumes that there is a single round of communication between the vertices (when each node collects the certificates of its neighbors), which is why this type of proof labelling scheme is called a \emph{one-round} proof labelling scheme. There is a more general definition \cite{pls} in which each node is allowed to gather the certificates of its neighbors at distance $t$, for some integer $t\ge 1$, but in our case it is enough to restrict ourselves to $t=1$.

\section{Cellular embeddings and rotation systems}\label{sec:surfaces}

\subsection{Surfaces}

\begin{figure}[htb]
 \centering
 \includegraphics[scale=1.2]{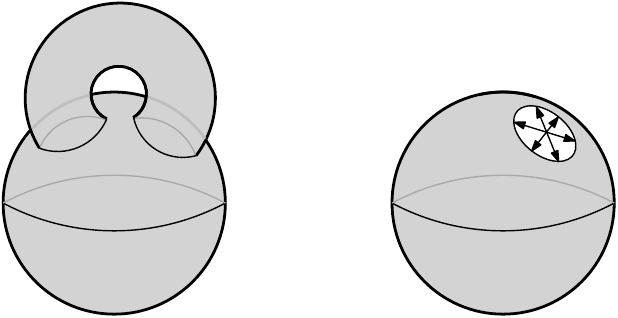}
 \caption{A sphere with a handle (left) and a sphere with a cross-cap (right).}
 \label{fig:handle}
\end{figure}

We refer the reader to the book by Mohar
and Thomassen~\cite{MoTh} for more details or any notion not defined
here. A {\em surface} is a non-null compact connected
2-manifold without boundary (meaning that every point has a neighborhood that is homeomorphic to an open subset of the plane). 
By the classification theorem of surfaces,
any surface is homeomorphic to one of the following  (see Figure~\ref{fig:handle}):

\begin{itemize}
    \item the \emph{orientable
  surface of genus~$h$}, obtained by adding $h\ge0$
\emph{handles} to the sphere.
    \item the \emph{non-orientable
  surface of genus~$k$}, 
    obtained by adding $k\ge1$
\emph{cross-caps} to the sphere (where a \emph{cross-cap} is a hole in the surface where antipodal points are identified).
\end{itemize}

\begin{figure}[htb]
 \centering
 \includegraphics[scale=0.8]{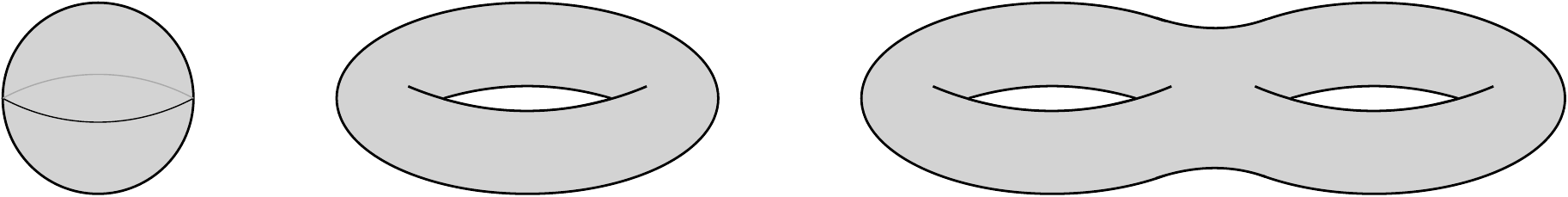}
 \caption{From left to right: the sphere, the torus, and the double-torus.}
 \label{fig:torus}
\end{figure}

The Euler genus of a surface $\Sigma$, denoted by $\mathbf{eg}(\Sigma)$, is equal to twice its genus if $\Sigma$ is orientable, and as its genus otherwise. So, as already mentioned
in the introduction the orientable surfaces of Euler genus 0, 2 and 4 are respectively the sphere (or plane), the torus and the double torus (see Figure~\ref{fig:torus}).

An \emph{embedding} of a graph $G$ on a surface $\Sigma$ is a representation of $G$ on $\Sigma$ where the vertices of $G$ are distinct elements of $\Sigma$, and every edge of $G$ is a simple arc connecting in $\Sigma$ the two vertices which it joins in $G$, with the property that the interior of every edge is disjoint from other edges and vertices. Removing the vertices and edges of $G$ from the surface $\Sigma$ divides $\Sigma$ into connected components, called \emph{faces} of the embedding of $G$ (or faces of $G$, if the embedding is clear from the context). We say that an embedding is \emph{cellular} if every face is
homeomorphic to an open disk of~$\mathbb{R}^2$. 
If $G$ is   a graph with a cellular embedding in a surface of Euler genus $g$, then Euler's formula states that \[ |V(G)|-|E(G)|+|F(G)|=2-g,
\] where $V(G)$, $E(G)$, and $F(G)$ denote the set of vertices, edges, and faces of (the embedding of) $G$.

\medskip

In order to prove Theorem~\ref{thm:main}, it will be convenient to restrict ourselves to cellular embeddings. We will actually prove the following.

\begin{theorem}\label{thm:cellular}
For any integer $g\ge 0$, the class of graphs with a cellular embedding on an orientable (resp.\ non-orientable) surface of Euler genus at most $g$
has a one-round proof labelling scheme with complexity at most $O(\sqrt{g}\cdot \log n)$.
\end{theorem}

\noindent \emph{Proof of Theorem~\ref{thm:main} assuming Theorem~\ref{thm:cellular}.} A graph is embeddable on an orientable surface $\Sigma$ of Euler genus $g$ if and only if it has a  cellular embedding in an orientable surface of Euler genus at most $g$ (see~\cite[Proposition 3.4.1]{MoTh}), so the orientable case of Theorem~\ref{thm:main} follows directly from the orientable case of Theorem~\ref{thm:cellular}. 

A graph $G$ is embeddable on a non-orientable surface $\Sigma$ of Euler genus $g$ if and only if it has a  cellular embedding in a non-orientable surface of Euler genus at most $g$ or $G$ is a tree (see~\cite[Proposition 3.4.2]{MoTh}). As trees have a simple one-round proof labelling scheme with complexity at most $O(\log n)$~\cite{AKY97}, the non-orientable case of Theorem~\ref{thm:main} can also be deduced directly from the non-orientable case of Theorem~\ref{thm:cellular}. 
\hfill $\Box$

\medskip

A graph $G$ is said to be \emph{$k$-degenerate} if there is an ordering $v_1,\ldots,v_n$ of the vertices of $G$, such that for any $1\le i \le n$, the vertex $v_i$ has at most $k$ neighbors $v_j$ with $j<i$. Note that the notion of degeneracy is very similar to the notion of \emph{edge-arboricity}, which has been extensively studied in the context of distributed algorithms (the two parameters are within a multiplicative constant of each other). 
Using Euler's formula, it is not difficult to derive the following result due to Heawood (see~\cite[Theorem 8.3.1]{MoTh}).

\begin{theorem}[Heawood]\label{thm:heawood}
For every $g\ge 0$, every graph embeddable on a surface of Euler genus at most $g$ is $k$-degenerate, with $k=\max\left(5,\tfrac{5+\sqrt{1+24g}}2\right)$.
\end{theorem}

In the next section we explain how the topological aspects of cellular embeddings can be translated into purely combinatorial notions in the case of orientable surfaces. The generalization to non-orientable surfaces is explained in Section~\ref{sec:no}.

\subsection{Rotation systems}

Let $G$ be a graph. A \emph{half-edge} of $G$ is a pair  $(v,e)$, where $v\in V(G)$ and $e$ is an edge incident to $v$. We say that $(v,e)$ and $v$ are incident. The set of all half-edges of $G$ is denoted by $B(G)$. A \emph{rotation system} of $G$ is a pair of permutations $(\sigma,\alpha)$ acting on $B(G)$, such that 
\begin{itemize}
    \item for any edge $e=uv\in E(G)$, $\alpha(v,e)=(u,e)$ and $\alpha(u,e)=(v,e)$ (i.e., $\alpha$ is an involution with no fixed point), and
    \item for each orbit of $\sigma$, there is a vertex $v\in V(G)$ such that the orbit consists of all the half-edges incident to $v$ (in other words, we can view $\sigma$ as a circular order on the half-edges incident to each vertex of $G$).
\end{itemize}

Each cellular embedding of a graph $G$ in some orientable surface $\Sigma$ can be translated into a rotation system by defining $\alpha$ as above and $\sigma$ as the collection of circular orders on the half-edges around each vertex, in the positive orientation of the surface. Note that each orbit of $\sigma \circ \alpha$ corresponds to a different face of the embedding (where the half-edges appear in the negative orientation of the surface).

\begin{figure}[htb]
 \centering
 \includegraphics[scale=1.2]{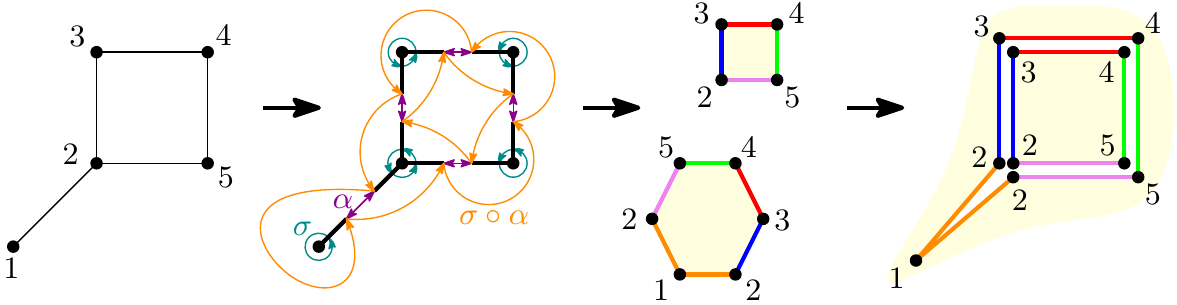}
 \caption{The conversion of a graph embedded in the plane (or the sphere) into a rotation system, and back into a graph embedded in the sphere (obtained by gluing the two polygons on edges of matching colors).}
 \label{fig:rotation}
\end{figure}

The Heffter-Edmonds-Ringel rotation principle (see Theorem 3.2.4 in~\cite{MoTh}) states that every cellular embedding of a connected graph $G$ in an orientable surface is uniquely determined, up to oriented homeomorphism, by its rotation system.

\medskip

Although it will not be needed in the remainder, it is worth explaining how to recover the embedding of $G$ in $\Sigma$ from the rotation system $(\sigma,\alpha)$ with ground set $B$. Each orbit $\sigma_v$ of $\sigma$ is associated to a distinct vertex $v$, and each orbit $\alpha_e$ of $\alpha$ is associated to a distinct edge $e$ connecting the two vertices associated to the two elements of $\alpha_e$. The graph resulting from this construction is precisely $G$. To each  orbit  of $\sigma \circ \alpha$, we associate a polygon $f$ whose sides are indexed by the edges of $G$  and whose vertices are indexed by the vertices of $G$. Note that vertices and edges of $G$ might appear several times on the same polygon or on different polygons (indeed, each edge appears twice and each vertex $v$ appears $d(v)$ times among all the polygons). The circular order on the vertices and edges on each polygon in the negative orientation coincides with the circular order of the elements of $B$ in the corresponding orbit of $\sigma \circ \alpha$.  For any edge $e$ of $G$, we glue the two polygons containing $e$ together on $e$ (if a single polygon contains $e$ twice, we glue the two sides corresponding to $e$ together), by respecting the natural orientation of $e$ (that is, if $e=uv$, the vertex $u$ of one polygon is identified with the vertex $u$ of the other polygon, and similarly for $v$), see Figure~\ref{fig:rotation}.

\section{Certificates for orientable surfaces}\label{sec:o}

%A graph is embeddable on some orientable surface $\Sigma$ if and only if it is embeddable on an orientable surface of Euler genus at most $g=\mathbf{eg}(\Sigma)$.
%We now describe our certificates for being embeddable on an orientable surface  of Euler genus at most $g$, as well as the verification part at each vertex. 

In this section we prove the orientable case of Theorem~\ref{thm:cellular}.
Recall that all graphs in this paper are assumed to be simple and connected.

%\textcolor{blue}{remarque que pour $g>0$ on doit ajouter "connexe" sinon ça n'est pas possible de faire ca localement, cf enonce du theoreme}

\medskip

Let $G$ be a graph with a cellular embedding on an orientable surface $\Sigma$ of Euler genus at most $g$.  Let $(\sigma,\alpha)$ be the rotation system associated to some cellular embedding of $G$ in $\Sigma$, and let $F(G)$ denote the set of faces of the embedding.

\medskip

Our certificate consists in two parts: (1) a (local) description of the rotation system $(\sigma,\alpha)$, and (2) a spanning tree of $G$ which supports the (local) computation of the Euler genus, via Euler's formula. Recall that according to Euler's formula, the value of $g$ can be deduced from the values of the number of vertices, edges, and faces of the embedding. In order to compute the number of faces we will also need to store (locally) a circular order on the edges of each face. This will allow us to choose a root edge or vertex on each face, and avoid multiple counting of the same face along our spanning tree.

\subsubsection*{Edge certificates} To describe our certificates it will be convenient to assume in the remainder that not only vertices, but also edges are given certificates, and that each vertex collects the certificates of its neighboring vertices and incident edges before choosing whether to accept or reject the instance. 

Assume that $G$ is a $k$-degenerate graph (with an ordering of its vertices witnessing this property), and each vertex $v$ and edge $e$ is assigned a certificate $c(v)$ or $c(e)$ of size at most $t$. Then each vertex $u$ can store the certificates of the (at most $k$) edges $uv$ such that $v$ lies before $u$ in the order. We are now in the case where only vertices store certificates (of size at most $t+k(t+\log n)$), and after collecting the new certificates of its neighbors, each vertex has access to the certificates $c(v)$ of its neighbors $v$ and $c(e)$ of its incident edges $e$. Consequently, it follows from Theorem~\ref{thm:heawood} that in order to prove Theorem~\ref{thm:main}, it is enough to assign certificates of size $O(\log n)$ to the vertices and edges of $G$ (and to assume that each vertex collects the certificates of its adjacent vertices and incident edges before accepting or rejecting the instance). The same observation was used previously in~\cite{planar,genus} (see also~\cite{BFT}).

\subsection{A distributed rotation system}

Let $v\in V(G)$. The certificate of $v$ contains the identifier $\mathrm{id}(v)$ of $v$. Since for any vertex $u$ in $G$, $\mathrm{id}(u)\in \{1,\ldots,\text{poly}(n)\}$, storing a constant number of identifiers takes $O(\log n)$ bits.

\medskip

In the remainder, it will be convenient to talk about the \emph{identifier} $\mathrm{id}((v,e))$ of a half-edge $(v,e)$, which we define as the the pair $(\mathrm{id}(v),\mathrm{id}(u))$, where $u$ is the endpoint of $e$ distinct from $v$.

\subsubsection*{Certificates around the vertices} 
For each vertex $v$, let $(v,e_0),(v,e_1),\ldots,\allowbreak (v,e_{d(v)-1})$ be the half-edges incident to $v$ in the positive orientation, starting with some arbitrary half-edge $(v,e_0)$ incident to $v$. For each $0\le i \le d(v)-1$, we say that $(v,e_i)$ has \emph{$v$-index $i$}, where $v$-indices are understood modulo $d(v)$, and we denote this half-edge by $\langle v\rangle_i$. By extension, we also say that the endpoint of $e_i$ distinct from $v$ has \emph{$v$-index $i$}, so that the circular ordering of the half-edges around $v$ coincides with a circular ordering of the neighbors of $v$ (this is possible here since we deal with simple graphs). Then each edge $uv$ is given as certificate the identifiers of $u$ and $v$, together with the $u$-index of $v$ and the $v$-index of $u$.

\subsubsection*{Certificates around the faces} 
Consider a half-edge $(v,e)$. The face $f$ associated to the orbit of $\sigma\circ \alpha$ containing $(v,e)$ is said to be the \emph{face bounding the half-edge $(v,e)$}, and we say that $(v,e)$ is bounded by $f$. The half-edge $(\sigma\circ\alpha)(v,e)$ is called \emph{the next half-edge on $f$} with respect to $(v,e)$. Note that if an edge $e$ of $G$ is incident to a single face $f$, the two half-edges of $e$ are bounded by $f$, while if $e$ is incident to two distinct faces $f_1,f_2$, one half-edge of $e$ is bounded by $f_1$ and the other is bounded by $f_2$.

\medskip

For each face $f$ of $G$, the prover  considers an arbitrary half-edge bounded by $f$ and sets it as the \emph{root} of $f$ (in the remainder, if the root of $f$ is $(v,e)$, we say that \emph{$f$ points to $v$}). The prover then
assigns integers to the half-edges bounded by $f$ as follows: for any half-edge $(v,e)$ bounded by $f$, the $f$-index of $(v,e)$ is the smallest integer $i\ge 0$ such that $(v,e)=(\sigma \circ \alpha)^{i}(u,e_0)$, where $(u,e_0)$ denotes the root half-edge of $f$. So the root half-edge of $f$ has $f$-index 0, and the maximum $f$-index is  $d(f)-1$,  where $d(f)$ denotes the degree of $f$ (the number of edges in a boundary walk of $f$, where edges appearing twice in the walk are counted with repetition). Note that if some half-edge has $f$-index $i$, the next half-edge on $f$ has $f$-index $i+1$ if and only if it is different from the root half-edge of $f$.

Now, consider any edge $uv$, and let $f_u$ be the face bounding the half-edge $(u,uv)$ and let $f_v$ be the face bounding the half-edge $(v,uv)$. Then the edge $uv$ is also given as certificate the $f_u$-index of  $(u,uv)$ together with the identifier of the root half-edge of $f_u$, and similarly the $f_v$-index of  $(v,uv)$ together with the identifier of the root half-edge of $f_v$.

\medskip

Each vertex has degree at most $n-1$ and each face has degree at most $2|E(G)|=O(n^2)$, so storing each $v$-index or $f$-index takes at most $O(\log n)$ bits. It follows that storing all the information described above  takes at most $O(\log n)$ bits per vertex and edge, and thus $O(k\log n)$ per vertex (where $k=\max\left(5,\tfrac{5+\sqrt{1+24g}}2\right)$, 
see Theorem~\ref{thm:heawood}).

\medskip

We now describe the verification process at each vertex.

\subsubsection*{Verifying the vertices} Each vertex collects the certificates of all its neighbors and incident edges. After having collected these certificates, both $u$ and $v$ are supposed to have all the \emph{information concerning the edge $uv$}, namely: the identifiers $\mathrm{id}(u)$ and $\mathrm{id}(v)$, the $u$-index $i$ of $v$, the $v$-index $j$ of $u$, the identifiers of the root half-edges of the faces $f$ and $f'$ bounding $\langle u\rangle_{i}$ and $\langle v\rangle_{j}$ respectively, the $f$-index of $\langle u\rangle_{i}$ and the $f'$-index of $\langle v\rangle_{j}$.
The verifier at each vertex $v$ checks that  the set of $v$-indices of the neighbors of $v$ forms a circular permutation of $\{0,\ldots,d(v)-1\}$, thus certifying that the information collected by $v$ is consistent with the local view of $v$ in some embedding of $G$. 

\medskip

\begin{figure}[htb]
 \centering
 \includegraphics[scale=1.2]{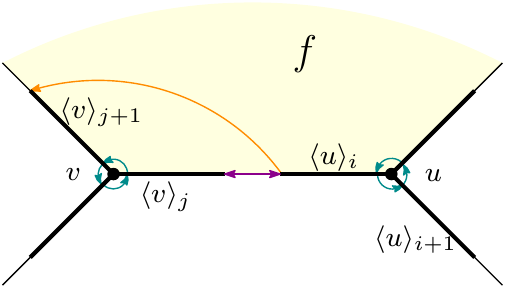}
 \caption{The next half-edge on $f$.}
 \label{fig:face}
\end{figure}

Let $(\sigma,\alpha)$ be the rotation system given by the $v$-indices of incident half-edges at each vertex $v$.
By the Heffter-Edmonds-Ringel rotation principle, $(\sigma,\alpha)$ defines a unique cellular embedding of $G$ on an orientable surface $\Sigma$ (up to oriented homeomorphism). Note that for any edge $uv$ where $i$ is the $u$-index of $v$ and $j$ is the $v$-index of $u$, and $f$ is the face bounding $\langle u\rangle_{i}$ in $\Sigma$, the next half-edge on $f$ with respect to $\langle u\rangle_{i}$ is $\langle v\rangle_{j+1}$ (see Figure~\ref{fig:face}).

\subsubsection*{Verifying the faces}
For each vertex $v$, and each $0\leq j \leq d(v)-1$ we do the following. Let $u$ be the neighbor of $v$ with $v$-index $j$, and let $i$ be the $v$-index of $u$ (see Figure~\ref{fig:face}). The verifier at $v$ checks that the half-edges $\langle u\rangle_{i}$ and $\langle v\rangle_{j+1}$ agree on the identifier of the root half-edge of the face $f$ bounding them, so that  the knowledge of the root half-edge of $f$ is consistent along the face $f$. Since this verification is done by all vertices on the boundary of $f$, if no vertex rejects the instance, then each face has a unique root half-edge.
    In order to make sure that this root half-edge of $f$ is actually  bounded by $f$, the verifier at  $v$ simply checks that the $f$-index of $\langle v\rangle_{j+1}$ is equal to $0$ if $\langle v\rangle_{j+1}$ is the root half-edge of $f$, or equal to $1$ plus the $f$-index of $\langle u\rangle_{i}$ otherwise. Since the face $f$ is finite and circularly ordered, some half-edge $(u,e)$ bounded by $f$ must have $f$-index at least the $f$-index of the next half-edge $(u',e')$ on $f$, and by definition this is only possible if $(u',e')$ has $f$-index 0. It follows that if the verifier at each vertex agrees with the instance, each face $f$ has  a unique root half-edge, and this root half-edge is bounded by $f$ (so $f$ points to a unique vertex, and this vertex is lying on $f$).
    
    Hence, if no vertex has rejected the instance so far, each face bounds a unique root half-edge, and each vertex knows whether each of its half-edges is the root half-edge of the face bounding it. In particular, each vertex knows the number of faces pointing to it.

\subsection{Computation of the Euler genus}
\label{sec:genus}

Using the information collected by each vertex $v$, and assuming all vertices have accepted the instance so far, we now certify that the surface $\Sigma$ has Euler genus at most $g$. To do this, it suffices to compute $|V(G)|$, $|E(G)|$, and $|F(G)|$ and apply Euler's formula. We will do so by collecting the number of vertices, edges and faces along a spanning tree. Let $T$ be a rooted spanning tree in $G$ with root $r$. This spanning tree is certified locally using the following classical scheme (see~\cite{AKY97,APV91,IL94}): the prover gives  the identifier $\mathrm{id}(r)$ of the root of $T$ to each vertex $v$ of $G$, as well as $d_T(v,r)$, its distance to $r$ in $T$, and each vertex $v$ distinct from the root is also given the identifier of its parent $p(v)$ in $T$. The verifier at $v$ starts by checking  that $v$ agrees with all its neighbors in $G$ with the identity of the root $r$ of $T$. If so, if $v\ne r$, $v$ checks that $d_T(v,r)=d_T(p(v),r)+1$.
Once the rooted spanning tree $T$ has been certified, each vertex of $G$ knows its children in $T$. This can be used to check that
\begin{itemize}
    \item all vertices agree on the same number $n=|V(G)|$ of vertices: In order to do this, the prover gives $n$ to each vertex $v$ of $G$, as well as a counter $\nu(v)$ which is equal to  the number of vertices in the subtree of $T$ rooted in $v$. The verifier at every vertex $v$ simply checks that $v$ has the same value of $n$ as its neighbors in $G$, and 
    that $\nu(v)$ is equal to 1 plus the sum of $\nu(u)$, for all children $u$ of $v$ (if any). Note that this can be checked locally. It only remains to check that for the root $r$ of $T$, $\nu(r)=n$.
    \item all vertices agree on the same number $m=|E(G)|$ of edges: Again, the prover gives the value of  $m$ to each vertex $v$ of $G$, together  with a counter $\mu(v)$ defined as the half of the sum of the degrees $d_G(u)$ of the vertices $u$ in the subtree of $T$ rooted in $v$. The verifier at $v$ only needs to check that it agrees on the value of $m$ with its neighbors in $G$, and that $\mu(v)$ is $\tfrac12 d_G(v)$ plus the sum of $\mu(u)$, for all children $u$ of $v$ (if any). Since
    $m=|E(G)|=\tfrac12 \sum_{v\in V(G)}d_G(v)$,
    it remains to check that for the root $r$ of $T$, $\mu(r)=m$. 
    \item all vertices agree on the same number $|F(G)|$ of faces: Again, the prover gives the value of $|F(G)|$ to each vertex $v$ of $G$, together  with a counter $\phi(v)$ equal to the number of faces pointing to vertices lying in the subtree of $T$ rooted in $v$ (recall that each face has a unique root half-edge, and each vertex knows the number of faces pointing to it). The verifier at $v$ checks that $\phi(v)$ is the number of faces pointing to $v$ plus the sum of $\phi(u)$, for all children $u$ of $v$ (if any). It remains to check that for the root $r$ of $T$, $\phi(r)=|F(G)|$. 
\end{itemize}

It follows that, assuming no vertex has rejected the instance so far,  each vertex has now access to $|V(G)|$, $|E(G)|$, and $|F(G)|$, and can check whether \[ 2+|E(G)|-|V(G)|-|F(G)|\le g.\] As a consequence of Euler's formula, this is equivalent to saying that the rotation system associated to $G$ embeds $G$ in an orientable surface of Euler genus at most $g$. 

\smallskip

This concludes the proof of Theorem~\ref{thm:main} for orientable surfaces.

\section{Non-orientable surfaces}\label{sec:no}

In this section we explain how to prove the non-orientable case of Theorem~\ref{thm:cellular}.
The case of non-orientable surfaces is very similar to the case of orientable surfaces, but there is an additional twist. An \emph{embedding scheme} is a rotation system $(\sigma,\alpha)$, except that each orbit $e$ of $\alpha$ has a sign $\lambda_e\in \{-1,1\}$. Given a cellular embedding of a graph $G$ in a surface $\Sigma$ (which is orientable or non-orientable), we can associate a circular order on the half-edges incident to each vertex $v$, by choosing an arbitrary orientation of the topological neighborhood of $v$ (positive or negative). This choice of local orders around the vertices gives $\sigma$, and the edges give $\alpha$, as before. Since we have chosen arbitrary orientations around the vertices, the orientations around two adjacent vertices $u$ and $v$ may not be consistent (i.e., agree on a small topological neighborhood around the edge $uv$). If they are consistent we set $\lambda_{uv}=1$ and otherwise we set $\lambda_{uv}=-1$. The surface is orientable if and only if there is a choice of local orientations that is globally consistent, that is such that the resulting signs satisfy $\lambda_e=1$ for every edge $e$. The surface is non-orientable if and only if $G$ contains a cycle $C$ which has odd number of edges $e$ with $\lambda_e=-1$ (see Section 3.3 in~\cite{MoTh}).

\begin{figure}[htb]
 \centering
 \includegraphics[scale=1.2]{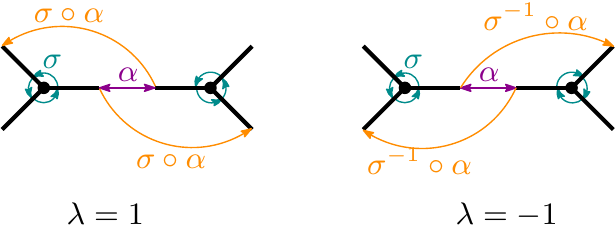}
 \caption{Description of the next half-edge on a face depending whether the local orderings of two adjacent vertices are consistent (left) or not (right).}
 \label{fig:negative}
\end{figure}

It turns out that the Heffter-Edmonds-Ringel rotation principle still holds in this more general setting (see Theorem 3.1.1 in~\cite{MoTh}). The only difference when retrieving the faces of the embedding is that in order to find the next half-edge on some face, with respect to some half-edge $(v,e)$, we consider the sign $\lambda_e$. If $\lambda_e=1$, the next half-edge on the face bounding $(v,e)$ is $(\sigma\circ \alpha)((v,e))$, as before.  If $\lambda_e=-1$, the next half-edge on the face bounding $(v,e)$ is $(\sigma^{-1}\circ \alpha)((v,e))$. So instead of identifying the faces in the embedding of $G$ with orbits of $\sigma\circ \alpha$ as before, we identify them with orbits of the function $\varphi: (v,e)\mapsto (\sigma^{\lambda_e}\circ \alpha)((v,e))$. This is illustrated in Figure~\ref{fig:negative}. 

\medskip

To adapt the certificate of the orientable setting to this more general framework,  the value of $\lambda_{uv}$ is added to the certificate of each edge $uv$. Using this additional information, the next half-edge on a face is computed using $\varphi$ instead of $\sigma\circ \alpha$. A single bit is added to the certificate of every edge, and so at most $k=O(\sqrt{g})$ bits are added to  the certificate of every vertex. 

\medskip

At this point, assuming no vertex has rejected the instance so far, we have certified that $G$ can be embedded in some surface of Euler genus at most $g$ (orientable or not). To conclude, it remains to certify that $G$ can be embedded in a \emph{non-orientable} surface of Euler genus at most $g$ (we thank an anonymous referee for pointing out that this final step was missing in an earlier version of this manuscript). As suggested above, this can be done by certifying the existence of a cycle $C$ which has an odd number of \emph{negative} edges, i.e.\ edges $e$ with $\lambda_e=-1$. To do so, the prover  chooses the rooted spanning tree $T$ of Section~\ref{sec:genus} such that 
\begin{itemize}
    \item its root $r$ is in a cycle $C$ with an odd number of negative edges,
    \item there is a negative edge $e_r$ of $C$ that is incident to $r$ and not in $T$, and
    \item the path $C\setminus e_r$ is a path in $T$.
\end{itemize}  
Then each vertex $v$ stores the parity of the number of negative edges on the path from $v$ to $r$ in $T$. This can be certified by similar techniques as in  Section~\ref{sec:genus}, namely with a counter $\eta$ at each vertex such that $\eta(r)=0$ for the root $r$, and, for each vertex $v\in V(G)\setminus r$, the value of $\eta(v)$ is equal to $\eta(p(v))+1$ if $\lambda_{vp(v)}=-1$ and $\eta(p(v))$ otherwise. Finally the verifier at $r$ checks that the value $\eta(v)$
of its neighbor $v$ along $e_r$ is equal to $0$. With the additional edge $e_r$, this gives an odd number of negative edges in $C$, ensuring that the surface is non-orientable.

At most $O(\log n)$ bits are added to  the certificates of the vertices to certify $C$, so the complexity of the proof labelling scheme remains $O(\sqrt{g}\log n)$. 

\section{Conclusion and open problems}\label{sec:ccl}

In this paper we gave a simple proof of Theorem~\ref{thm:main}, which shows that embeddability in a fixed surface (and in particular planarity) has a one-round proof labelling scheme with complexity $O(\log n)$, which is best possible. Graphs embeddable on a fixed surface form important examples of \emph{minor-closed classes}, that is classes $\mathcal{F}$ such that for any graph $G\in \mathcal{F}$, any minor of $G$ (i.e., any graph obtained from $G$ by deleting vertices and edges and contracting edges) lies in $\mathcal{F}$.

The authors of~\cite{planar,genus} asked whether any minor-closed class $\mathcal{F}$ has a one-round proof labelling scheme with complexity $O(\log n)$, and noted that even giving local certificates of $O(n^\alpha)$ bits for such classes, for some $\alpha<1$, seems to be a difficult challenge. It can be observed that on the other hand, it follows from the Graph Minor Theorem of Robertson and Seymour~\cite{RS04} that for any proper minor-closed class $\mathcal{F}$, there is a one-round proof labelling scheme with complexity $O(\log n)$ for the \emph{complement} of $\mathcal{F}$. This is because such a class $\mathcal{F}$ has a finite number of minimal obstructions, and these obstructions (if they appear) can be easily certified using $O(\log n)$ bits per vertex (see~\cite{planar}, where constructions are given for specific minors). Note that compact proof labelling schemes for $H$-minor free graphs (when $|V(H)|\le 4$) were given in~\cite{BFT}.

We note here that even proving that graphs of bounded treewidth have one-round proof labelling schemes with complexity $O(\log n)$ seems to be non-trivial. 

\begin{acknowledgement}
We thank the reviewers for their comments and suggestions.
\end{acknowledgement}

\end{document}